\documentclass[]{IEEEtran}
\usepackage{xcolor}
\usepackage{soul}
\usepackage{booktabs,siunitx}
\usepackage[cmex10]{amsmath}
\usepackage{amsfonts}
\usepackage{graphicx,caption,subcaption} 
\usepackage{sidecap}
\usepackage{mathtools}
\usepackage{enumerate}
\DeclareMathOperator{\Tr}{Tr} 
\usepackage{amsmath}%
\usepackage{MnSymbol}%
\usepackage{wasysym}%
\newtheorem{theorem}{Theorem} 
\newtheorem{lemma}{Lemma} 
\newtheorem{definition}{Definition} 

\usepackage{soul}
\usepackage{cite}

\definecolor{lightgray}{gray}{0.9}

\begin{document}
	\title{A Rapid Fault Reconstruction Strategy Using   a Bank of  Sliding Mode Observers}
	
	\author{
		  Mehran Shakarami,  Kasra~Esfandiari, Amir Aboulfazl Suratgar, and  Heidar Ali Talebi   
		\thanks{M. Shakarami  {is with the Engineering and Technology Institute Groningen,
				University of Groningen, 9747 AG Groningen, The Netherlands (e-mail:
				m.shakarami@rug.nl).}} 
				\thanks{K. Esfandiari is   with the Center for Systems Science, Yale University, New Haven, CT,  USA e-mail: (kasra.esfandiari@yale.edu).}
			\thanks{A.A. Suratgar and H.A. Talebi  {are with the Center of Excellence on Control and Robotics, Department of Electrical Engineering}, Tehran Polytechnic, Tehran, Iran (e-mail:
					\{a.suratgar;alit\}@aut.ac.ir).}} 
		 

	\markboth{ }%
	{Esfandiari \MakeLowercase{\textit{et al.}}: Bare Demo of IEEEtran.cls for IEEE Journals}

	\maketitle
	
	\begin{abstract}
	  This paper deals with the design of a model-based rapid fault detection and isolation strategy using sliding mode observers. To address this problem, a new    scheme is proposed by adaptively combining the information provided by a bank of   observers. In this regard, a new structure for sliding mode observers is considered. Then, the well-known recursive least square algorithm is utilized to merge individual state estimations suitably such that the system fault is detected faster. The required condition for enhancing   perfect state estimation is derived, and the stability of the overall system is proven via Lyapunov's direct method.   The supremacy of  proposed scheme is fully discussed  through mathematical analyses as well as simulations. 
		%
		
	\end{abstract}

\section{Introduction}

Reconstructing   unknown inputs of a given process is of great significance from both practical and theoretical aspects.  One example of unknown inputs is exogenous disturbances. It is well-known that if such disturbances are not identified and compensated properly, they will deteriorate performance of the closed-loop system. To tackle this problem, several disturbance observers have been developed in the control literature 
\cite{5., 1., autom4.}. System components fault, which frequently occurs in engineering systems, is another example of unknown inputs. Occurrence of fault may cause  irrecoverable damages and even failure of the whole process. To avoid such circumstances, one solution   is to consider  redundancy for  critical components of the process and always check the signals provided by these duplicated components. By using such structures and employing a voting mechanism, it is possible to detect faulty signals, and in turn,  select the appropriate signals. This strategy has been  widely used in industries to obtain a process with  high availability. Although such a simple approach has proven to provide reliable schemes in practice, one cannot turn a blind eye on the fact that it is costly and not energy efficient.   Moreover, in some applications it is not possible to duplicate  system components due to the nature of the understudy problem, space limitations,   accumulation of noise, etc.  To avoid such problems, in the past few decades, several attempts have been made to reconstruct faults using automatic strategies  \cite{talebi, 20., 23., robust2}. 
In \cite{20.}, Sliding Mode Observers (SMOs)  are employed to design a model-based structure to reconstruct    waste-gate  faults in  turbocharged gasoline engines.  In \cite{23.}, an algorithm  is presented for fault detection in linear time-varying discrete systems  using Luenberger  observers. Partial kernel principle component analysis is employed for health monitoring of aero-derivative industrial turbines in   \cite{44.}.  

Many existing   methodologies in control literature  assume  that all system states are accessible  \cite{NN}. However,   this assumption is not always realistic.   It is well-known that this  issue can be addressed by employing suitable observers. In the past decades, several observers and observer-based strategies have been  presented for linear and nonlinear systems  \cite{david1971introduction, edwards1998sliding,icee,  maghaleman}.        Among these observers, SMOs,  which  are capable of rejecting impacts of disturbances on the final estimations, have obtained great deal of attentions from  fault detection community.  This feature, disturbance rejection,  is  employed  for   fault detection, isolation, and input reconstruction in 
\cite{robust4, edwards2000sliding, robust1, MMSliding}. 

On the other hand, it should be   emphasized that   not only does late detection of such malfunctions    cause fault propagation and failure of the whole process, but  it also  put the process operators' lives in jeopardy. Hence, early detection of such malfunctions  has received a considerable attention, and several attempts have been made to design rapid fault detection methodologies  
 \cite{130., 111., 120.}.  In \cite{130.}, an intelligent  modular method is proposed for fast fault detection and classification in power systems.
A robust fast fault detection approach is presented for T-S fuzzy systems in \cite{111.}.   Assuming the time derivative of  system output,  $ \dot y(t)$, is available for measurement;  an algorithm is presented for rapid  actuator fault detection in \cite{120.}.  
On the other hand, it is well-known that by using multiple observers one can estimate  
system states with better transient response   \cite{  postoyan2015multi, chong2015parameter}, which can be employed for fault detection  
purposes. In \cite{han2012new}, a new identification scheme is presented for Linear Time Invariant (LTI) systems in order to 
improve the parameter estimation performance. In this regard,  
a convex combination of all information provided by multiple models is  
utilized to  estimate system unknown parameters. Applying this idea   to different  systems with unknown parameters (e.g., LTI systems \cite{chen2014combination}, linear systems with unknown periodic parameters \cite{Report2, NaderCDC},  nonlinear systems 
\cite{pandey2015adaptive, CDC}, and twin-rotor system \cite{pandey2016control})   has shown that this approach is also capable of providing satisfactory performance in various cases. 


In this paper, a novel fault reconstruction  scheme is presented using a convex combination of the state estimations obtained from a bank  of SMOs. The method is composed of a new structure for multiple SMOs and the Recursive Least Squares (RLS) algorithm. The main  contributions of this paper are as follows: I) A new structure for state estimation of systems with unknown inputs is presented by adaptively combining multiple SMOs state estimations. 
II) Using the properties of convex combination, the existence of some unknown constant parameters  that provide a perfect state estimation is guaranteed. III) The stability of the  proposed scheme, consisting of interconnection of multiple dynamical systems,   is investigated, and it is proved that  the estimations converge to the actual values. IV) The mathematical performance analysis of the proposed scheme is provided which shows that the presented strategy is capable of resulting in  a better transient response in comparison to the conventional SMOs.  V) The presented fault reconstruction methodology does not require   accessibility of all system states.  

The remainder of paper is structured as  follows: The problem statement and SMO are presented in Section 2. In Section 3, the proposed observer scheme is introduced; then, the existence of a perfect state estimation is assured. Furthermore,  the stability of proposed structure and its performance  are also  discussed in this section. Simulation results are included in Section 4, and finally, Section 5 provides the conclusions.            

\section{Sliding Mode Observers} 
Consider the following uncertain system
\begin{equation}
\begin{split}
\dot{x}&=Ax+Bu+D\xi(t,x,u),\\
y&=Cx,
\end{split}
\label{sys}
\end{equation}
where $ x\in {\mathcal{R}}^n $ is the state vector, $ u\in \mathcal{R}^m $ is the bounded control input  which stabilizes the system, $ y\in \mathcal{R}^p $ is the system output, and $\xi(t,x,u)\in \mathcal{R}^q$ denotes an unknown bounded function satisfying the following inequality
\begin{equation}\label{up-xi}
\|\xi(t,x,u)\|\leq \overline{\xi},
\end{equation}
Moreover, it is assumed that $ p\geq q $ and  $ B $, $ C $, and $ D $ are full rank matrices \cite{edwards1998sliding}. 

To estimate the state of the system, a SMO with the following structure can be utilized 
\begin{equation}
\begin{split}
\dot{\hat{x}}&=A\hat{x}+Bu-G_l\tilde{y}+G_n\nu(\hat{x}),\\
\hat{y}&=C\hat{x},
\end{split}
\label{smo}
\end{equation}
where $\hat x$ and $\hat y $ represent estimations of $x$ and $y$, respectively,   $\tilde{y}=\hat{y}-y $ is the output estimation error, and $ \nu $ is a discontinuous term about the hyperplane
\begin{equation*}
\mathcal S_o=\{\tilde{x}\in\mathcal{R}^n:C\tilde{x}=0\}; \tilde{x}=\hat{x}-x ,
\end{equation*}
Moreover,  $ G_l \in \mathcal{R}^ {n\times p}$ and $ G_n \in \mathcal{R}^ {n\times p}$ are  gain matrices to be determined. Note that  in order  to estimate the system states precisely, it is  needed to reject the effects of unknown term $\xi $ on $\hat x$. 
The following definitions, lemmas, and  theorem are used throughout the paper. 

\begin{definition}[{\cite{bartolini2008modern}}]
	The Rosenbrock matrix $ \mathbb{R}(s) $ of the system $ (A,D,C) $  is given by
	\begin{equation}\label{rosen}
	\mathbb{R}(s)=\begin{bmatrix}
	sI_n-A & D\\
	C & 0
	\end{bmatrix}
	\end{equation}
	The values of $ s_0 $ such that $ \text{rank}(\mathbb{R}(s_0))<n+q $ are called invariant zeros of the system $ (A,D,C) $.
\end{definition}

\begin{definition}[\cite{bakelman2012convex}]
	A set $ \mathcal{K} $ in a linear space $ \mathcal{L} $ is called	convex if the line segment $ ab $ is contained in $ \mathcal{K} $ for any	elements $ a, b \in \mathcal{K}$, i.e., $ (1-\lambda)a+\lambda b\in \mathcal{K} $ for any pair $ (a,b)\in \mathcal{K} $ and any $ \lambda\in [0,1] $.
\end{definition}

\begin{lemma} [\cite{bakelman2012convex}]
	\label{lem:con-hull}
	Let $a_1,a_2,\cdots,a_m\in \mathcal{L}$ where $\mathcal{L}$ is a linear space. The intersection of all convex sets in $\mathcal{L}$ containing $a_i$s is called  the convex hull $\mathcal{K}$ of \{$a_i(i=1,2,\cdots,m)$\} and any element of which, $a'$,  can be expressed as follows:
	\begin{equation*}
	a'=\sum\limits_{i=1}^m\beta_i a_i
	\end{equation*}
	where $\beta_i \in [0,1]$ is a constant term satisfying $\sum\nolimits_{i=1}^{m}\beta_i =1$.
\end{lemma}

\begin{lemma} [\cite{hansen2015functional}]
	\label{lem:neum}
	Let $ F $ be a bounded linear operator on an arbitrary Banach space $ \mathcal{L} $. If $ \|F\|<1 $, then $ I-F $ has a bounded inverse as follows
	\begin{equation*}
	(I-F)^{-1}=I+\sum_{k=1}^{\infty}F^k
	\end{equation*}
\end{lemma}

\begin{theorem}[\cite{edwards1994development}]
	\label{th:smo-cond}
	Consider the uncertain system \eqref{sys}, SMO \eqref{smo}, and let
	\begin{equation*}
	\nu(\hat{x})=\begin{cases}
	-\rho\frac{P_2\tilde{y}}{\|P_2\tilde{y}\|} & \text{if }\tilde{y}\neq 0\\
	0 & \text{otherwise}
	\end{cases}	\end{equation*}
	where $ P_2\in \mathcal{R}^{ p\times p} $  and $ \rho\geq \overline{\xi}+\gamma_0 $ with $ \gamma_0>0 $. Then there exist matrices $ G_l $, $ G_n $, and $ P_2 $ such that $ A_0=A-G_lC $ is a Hurwitz matrix and the state estimation obtained from \eqref{smo}  converges to the state of the plant  if and only if:
	\begin{itemize}
		\item $ \text{rank}(CD)=q $.
		\item invariant zeros of the system $ (A, D, C) $ lie in the open left-hand side of complex plane, i.e., $\text{rank}(\mathbb{R}(s))=n+q$	for all $ s $ with non-negative real part.
	\end{itemize}
\end{theorem}

{In the sequel,  it is assumed that the conditions of  Theorem~\ref{th:smo-cond} are satisfied, and  a new observation scheme with more preferable performance is  presented. 
}
%
%
%
%
\section{The Main Results}
{In this section, the structure of the proposed observation scheme and its performance are fully discussed. Then, by utilizing the presented observation strategy,   the unknown input   reconstruction scheme as well as its stability analysis are addressed. }
\subsection{Structure of the Proposed Multiple Sliding Mode Observer}
In order to estimate the system states precisely with better transient response,   the estimated states are  provided by combining the observations obtained from multiple SMOs. 

%
%
%

The proposed observer scheme is composed of $N$ SMOs with the following dynamics, 
\begin{equation}
\begin{split}
\dot{\hat{x}}_i(\alpha,t)&=A\hat{x}_i(\alpha,t)+Bu(t)-G_l\tilde{y}_i(\alpha,t) \\&+G_n\nu(\sum\limits_{i=1}^{N}\alpha_i \hat{x}_i(\alpha,t),t),\\
\hat{y}_i(\alpha,t)&=C\hat{x}_i(\alpha,t),
\end{split}
\label{mult-smo}
\end{equation}
where $i=1, 2, \cdots, N$  ($N\geq n+1$), $ \tilde{y}_i(\alpha)=\hat{y}_i(\alpha)-y $ , $ \hat{x}_i(\alpha,0)=\hat{x}_i(0) $, and
\begin{equation*}
\nu(\sum\limits_{i=1}^{N}\alpha_i \hat{x}_i(\alpha))=\begin{cases} 
-\rho\frac{P_2 \sum\limits_{i=1}^{N}\alpha_i \tilde{y}_i(\alpha)}{\|P_2\sum\limits_{i=1}^{N}\alpha_i \tilde{y}_i(\alpha)\|} & \text{if }\sum\limits_{i=1}^{N}\alpha_i \tilde{y}_i(\alpha)\neq 0,\\
0 & \text{otherwise},
\end{cases}	\end{equation*}
with $0\le \alpha_i \le 1$, $\sum\nolimits_{i=1}^{N} {\alpha_i}=1$, and $\alpha=\begin{bmatrix}
\alpha_1 & \alpha_2 & \cdots & \alpha_{N}
\end{bmatrix}^T$. To use the full information provided by the aforementioned $N$ observers and obtain promising  estimation, the final state estimation $\hat x_o$  is considered as a convex combination of the above  estimations, i.e., 
\begin{equation}
\hat{x}_o(t)=\sum\limits_{i=1}^{N}\alpha_i \hat{x}_i(\alpha,t).
\label{xo}
\end{equation}

In the sequel, the analysis is divided into two parts: the algebraic part and the analytic part. In the former part, it is guaranteed that there exist  unknown fixed $\alpha_i^*$s such that by choosing $ \alpha_i=\alpha_i^* $s the perfect state estimation is obtained.  The latter part covers derivation of  appropriate adaptive laws  for estimating  these parameters.
The following lemma is presented to obtain the required condition for the existence of $ \alpha_i^* $s.

\begin{lemma}
	\label{lem:con-IC}
	If the initial conditions of $N$ ($N\geq n+1$) SMOs with the structure of  \eqref{mult-smo}, $ \hat{x}_i(0) $s, are chosen such that the initial condition of the uncertain system \eqref{sys}, $ x(0) $,  lies in their convex hull $ \mathcal{K}$, then there exist some   fixed $ \alpha_i^* $s such that
	\begin{equation}
	x(t)=\sum\limits_{i=1}^{N}\alpha_i^* \hat{x}_i(\alpha^*,t)
	\label{perfect}
	\end{equation}
	where $0\le \alpha_i^* \le 1$, $\sum\nolimits_{i=1}^{N} {\alpha_i^*}=1$, and $\alpha^*=\begin{bmatrix}
	\alpha_1^* & \alpha_2^* & \cdots & \alpha_{N}^*
	\end{bmatrix}^T$.
\end{lemma}

\textbf{Proof.} 
There exists a change of coordinates, using a non-singular  matrix $ J $  that transforms the triple $ (A,D,C) $ of system \eqref{sys} to new coordinates $ (\mathcal{A},\mathcal{D},\mathcal{C}) $ such that \cite{edwards1994development}: 
\begin{equation}\label{sys-trans}
\begin{split}
\dot{x}_I&=\mathcal{A}_{11}x_I+\mathcal{A}_{12}y+\mathcal B_1u\\
\dot{y}&=\mathcal A_{21}x_I+\mathcal A_{22}y+\mathcal B_2u+\mathcal D_2\xi
\end{split}
\end{equation}
where $ x_I\in \mathcal{R}^{(n-p)} $, $ y\in\mathcal{R}^p $,  $ \mathcal{A}=JAJ^{-1}$, $\mathcal{B}=JB$,  $\mathcal{C}=CJ^{-1} $, $\mathcal{D}=JD$, and $ \mathcal{A}_{11} $ is a stable matrix. 
On the other hand by using \eqref{mult-smo}, \eqref{xo}, and the fact that $ \alpha_i $s are constant, one can obtain the dynamics of state estimation as 
\begin{equation}
\begin{split}
\dot{\hat{x}}_o&=A\hat{x}_o+Bu-G_l\tilde{y}_o+G_n\nu(\hat{x}_o),\\
\hat{y}_o&=C\hat{x}_o,
\end{split}
\label{dynm-xo}
\end{equation}
where $\tilde{y}_o=\hat{y}_o-y$, $ \hat{x}_o(0)=\sum_{i=1}^{N}\alpha_i\hat{x}_i(0) $, and
\begin{equation*}
\nu(\hat{x}_o)=\begin{cases} 
-\rho\frac{P_2 \tilde{y}_o}{\|P_2 \tilde{y}_o\|} & \text{if } \tilde{y}_o\neq 0,\\
0 & \text{otherwise}.
\end{cases}	\end{equation*} 
Similarly, the following dynamics can be obtained from \eqref{dynm-xo} using the aforementioned transition matrix  $J$: 
\begin{equation}\label{smo-trans}
\begin{split}
\dot{\hat x}_I&=\mathcal{A}_{11}\hat x_I+\mathcal{A}_{12}\hat y+\mathcal B_1u(t)-\mathcal{G}_{l1}\tilde{y}_o\\
\dot{\hat y}&=\mathcal A_{21}\hat x_I+\mathcal A_{22}\hat y+\mathcal B_2u(t)-\mathcal{G}_{l2}\tilde{y}_o+\mathcal{G}_{n2}\nu
\end{split}
\end{equation}
where $ \mathcal{G}_l=\begin{bmatrix}
\mathcal{G}_{l1}^T & \mathcal{G}_{l2}^T
\end{bmatrix}^T =JG_l$ and $ \mathcal{G}_n=\begin{bmatrix}
0 & \mathcal{G}_{n2}^T
\end{bmatrix}^T=JG_n $. Let us select the design parameters as 
\begin{equation}\label{gain-part}
\begin{split}
\mathcal{G}_{l1}=\mathcal{A}_{12},\quad
\mathcal{G}_{l2}= \mathcal{A}_{22}-\mathcal{A}_{22}^s,\quad
\mathcal{G}_{n2}=\|\mathcal{D}_2\|I_p
\end{split}.
\end{equation}
where $\mathcal{A}_{22}^s$ is a design Hurwitz matrix.  Now, a  Lyapunov function candidate as $ V(\tilde{x}_I,\tilde{y}_o)=\tilde{x}_I^TP_1\tilde{x}_I+\tilde{y}_o^TP_2\tilde{y}_o $ with $ \tilde{x}_I=\hat{x}_I-x_I $ can be considered. In addition, two symmetric positive definite design matrices are defined as $ Q_1\in \mathcal{R}^{(n-p)\times (n-p)} $ and $ Q_2\in \mathcal{R}^{p\times p} $, and $ P_1 $ and $ P_2 $ are the  symmetric positive definite
matrices that satisfy the following Lyapunov equations
\begin{align}
P_2\mathcal{A}_{22}^s+(\mathcal{A}_{22}^s)^TP_2&=-Q_2\label{lyp-eq}\\
P_1\mathcal{A}_{11}+\mathcal{A}_{11}^TP_1&=-(\mathcal{A}_{21}^TP_2Q_2^{-1}P_2\mathcal{A}_{21}+Q_1)\nonumber
\end{align}
Then, following a procedure similar to what presented in \cite{edwards1994development}, it can be shown that $ \dot{V}<0 $, and in turn, one can  conclude that
\begin{equation}
V(t)<V(0)
\label{vt}
\end{equation}
To show that there exist some $\alpha_i^*$s satisfying \eqref{perfect}, 
one can  employ the transformation matrix $ J $ and rewrite the Lyapunov function $ V $ as follows:
\begin{equation}
\begin{split}
V&=\tilde{x}_o^TP\tilde{x}_o,\\
P&=J^T\begin{bmatrix}
P_1 & 0\\
0 & P_2
\end{bmatrix}J
\end{split}
\label{v-def}
\end{equation}
where $\tilde x_o=\hat x_o - x$ is the observation error. By  considering the  equation of 
\begin{equation}
\underline\lambda(P) \tilde{x}_o^T\tilde{x}_o\leq V\leq \overline\lambda (P) \tilde{x}_o^T\tilde{x}_o,
\end{equation}
where  $ \underline{\lambda}(P) $ and $ \overline{\lambda}(P) $ respectively are the smallest and largest eigenvalues of $ P $, one can rewrite \eqref{vt} as
\begin{equation}
\| \tilde{x}_o(t)\|\leq\sqrt{\frac{\overline\lambda(P)}{\underline\lambda(P)}}  \| \tilde{x}_o(0)\|.
\end{equation}
It is obvious that if $ \tilde{x}_o(0)=0 $, the right-hand side of the preceding inequality is zero, and consequently, $ \tilde x_o(t)$ will maintain zero  for all $ t\ge0 $.  Therefore, to conclude the proof, it is sufficient to  find the conditions under which $\tilde x_o(t=0)$ is zero. 
Towards this end, one can use  Lemma~\ref{lem:con-hull} and conclude that if the initial conditions of   SMOs, i.e., $\hat{x}_i(0)$s, are chosen such that $x(0)$ is in the convex hull $\mathcal{K}$ of $\hat{x}_i(0)$s, then  fixed $\alpha_i^*$s exist such that
\begin{equation*}
x(0)=\sum\limits_{i=1}^{N}\alpha_i^* \hat{x}_i(0)
\end{equation*}
Hence, one can compare the preceding equation with \eqref{xo} and see  that the required condition for the perfect state estimation,  $\tilde x_o(t=0)=0$, is satisfied. 
\hfill$ \square $

It is worth noting that for $ x(0) $ lies in the convex hull $ \mathcal{K} $ of $ \hat{x}_i(0) $s, one should at least utilize $ N=n+1 $ SMOs. 

So far, it has been  shown that there exist some $\alpha_i^*$s that provide perfect state estimation. However,  these parameters are required to be estimated since they are unknown. 
Let us define $ \tilde{x}_i=\hat{x}_i-x $, and use  \eqref{perfect} and $\sum\nolimits_{i=1}^{N}\alpha_i^* =1$ to obtain $ \sum\nolimits_{i=1}^{N}\alpha_i^* \tilde{x}_i(\alpha^*,t)=0 $. On the other hand, because $\alpha_N^*=1-\sum\nolimits_{i=1}^{N-1}\alpha^*_i$ and $ \tilde{x}_i(\alpha^*)-\tilde{x}_N(\alpha^*)=\hat{x}_i(\alpha^*)-\hat{x}_N(\alpha^*) $, one can get
\begin{equation}
\sum\limits_{i=1}^{N-1}\alpha_i^* (\hat{x}_i(\alpha^*,t)-\hat{x}_N(\alpha^*,t))=-\tilde{x}_N(\alpha^*,t)
\label{alpha-error-n}
\end{equation}
Now, one can employ \eqref{mult-smo} to obtain
\begin{equation}\label{e-xhat-dot}
\dot{\hat{x}}_i(\alpha^*,t)-\dot{\hat{x}}_N(\alpha^*,t)=(A-G_l C)(\hat{x}_i(\alpha^*,t)-\hat{x}_N(\alpha^*,t))
\end{equation}
It can be seen that $ \hat{x}_i(\alpha^*,t)-\hat{x}_N(\alpha^*,t) $ is obtained from a linear system  and in turn,  only depends on $ \hat{x}_i(0)-\hat{x}_N(0) $; therefore, if a matrix $ E $ is considered such that its $ i $th column is equal to $ \hat{x}_i(\alpha^*,t)-\hat{x}_N(\alpha^*,t) $, it is independent of $ \alpha^* $. As a result, \eqref{alpha-error-n} is rewritten as follows
\begin{equation}\label{reg-0}
E(t)\bar\alpha^*=-\tilde{x}_N(\bar\alpha^*,t)
\end{equation}
where  $\bar\alpha^*=\begin{bmatrix}
\alpha_1^* & \alpha_2^* & \cdots & \alpha_{N-1}^*
\end{bmatrix}^T$. If $ \tilde{x}_N(\alpha^*) $ was known, the previous equation could be utilized for estimating $ \bar\alpha^* $; however, it is unknown since $ x $ and $ \hat{x}_N(\bar\alpha^*) $ are not available. Nonetheless, one can premultiply \eqref{reg-0} by $ C $ and  get 
\begin{equation}\label{reg-1}
CE(t)\bar\alpha^*=-\tilde{y}_N(\bar\alpha^*,t)
\end{equation}
Since $ \hat{y}_N(\bar\alpha^*) $ is unknown in the the preceding equation, the RLS algorithm cannot be employed for estimating $ \bar\alpha^* $. However, it is proposed to employ  a modification of the RLS algorithm as follows
\begin{equation}
\begin{aligned}
\dot{\hat{\bar\alpha}}(t)&=-R(t) E(t)^TC^T(CE(t)\hat{\bar\alpha}(t)+\tilde{y}_N(\hat{\bar\alpha}(t),t)),\\
\dot{R}(t)&=-R(t)E(t)^TC^TCE(t)R(t),
\end{aligned}
\label{RLS}
\end{equation}
where $\hat{\bar\alpha}(0)=\hat{\bar\alpha}_0$, $ R(0)=\mu I$, $ \hat{\bar\alpha} $ is an estimation of $ \bar\alpha^* $, $ \tilde{y}_N(\hat{\bar\alpha})=\hat{y}_N(\hat{\bar\alpha})-y $, $ I $ is the identity matrix, and $ \mu $ is a positive constant. Moreover, since $ E(t) $ is independent of $\bar\alpha^* $, its $ i $th column is considered as  $ \hat{x}_i(\hat{\bar\alpha}(t),t)-\hat{x}_N(\hat{\bar\alpha}(t),t) $. Even though $ \tilde{y}_N(\hat{\bar\alpha}) $ is employed instead of $ \tilde{y}_N({\bar\alpha}^*) $, it will be shown later that the previous equation is able to provide an estimation of $ \bar\alpha^* $ suitable for estimating the state of the plant. 
Now, the obtained $ \hat{\bar\alpha} $ can be  exploited in the following SMOs for estimating the state variables
\begin{equation}
\begin{split}
\dot{\hat{x}}_i(\hat{\bar\alpha})&=A\hat{x}_i(\hat{\bar\alpha})+Bu-G_l\tilde{y}_i(\hat{\bar\alpha})+G_n\nu(\hat x_o),\\
\hat x_o&=\sum\limits_{i=1}^{N-1} \hat{\alpha}_i \hat{x}_i(\hat{\bar\alpha}) + (1-\sum\limits_{i=1}^{N-1} \hat{\alpha}_i) \hat{x}_N(\hat{\bar\alpha})
\end{split}
\label{last_msmo}
\end{equation}
where $\tilde{y}_i(\hat{\bar\alpha})= C\hat{x}_i(\hat{\bar\alpha})-y $ and
\begin{equation*}
\nu(\hat x_o)=\begin{cases} 
-\rho\frac{P_2 \sum\limits_{i=1}^{N}\hat\alpha_i \tilde{y}_i(\hat{\bar\alpha})}{\|P_2\sum\limits_{i=1}^{N}\hat\alpha_i \tilde{y}_i(\hat{\bar\alpha})\|} & \text{if }\sum\limits_{i=1}^{N}\hat\alpha_i \tilde{y}_i(\hat{\bar\alpha})\neq 0,\\
0 & \text{otherwise},
\end{cases}	\end{equation*}

Now, it is required to investigate the stability of the proposed observation scheme, constructed from two interconnected systems  \eqref{RLS} and \eqref{last_msmo}, as well as the estimation error. In this regard, the following theorem is presented.

\begin{theorem}
	For system \eqref{sys}, consider  the modified RLS algorithm \eqref{RLS} and $ N $ SMOs \eqref{last_msmo}. If a non-singular matrix $ J $ that transforms \eqref{sys} to \eqref{sys-trans} is considered,  the design matrices are chosen as
	\begin{equation}\label{gl-gn}
	G_l=J^{-1}\begin{bmatrix}
	\mathcal{A}_{12}\\ \mathcal{A}_{22}-\mathcal{A}_{22}^s
	\end{bmatrix},G_n=\|\mathcal{D}_2\|J^{-1}\begin{bmatrix}
	0\\  I_p
	\end{bmatrix}
	\end{equation}
	and $ P_2 $ satisfies \eqref{lyp-eq}, then it can be guaranteed that the observation error system is quadratically stable. Moreover, $ \hat{x}_i $s are uniformly ultimately bounded, and $\hat{\bar\alpha}$ and $R$ are bounded.
	\label{th:conver}
\end{theorem}

\textbf{Proof.} Throughout the proof, it is required to employ the upper bounds of $ \|E(t)\| $ and $ \|R(t)\| $. In this regard, one can employ the definition of $ E $ and \eqref{e-xhat-dot} to obtain
\begin{equation}\label{e-dot}
\dot{E}(t)=A_0 E(t)
\end{equation}
where $ A_0=A-G_l C $ is a Hurwitz matrix. Therefore, we have $ E(t)=e^{A_0 t}E(0) $, and $ \|E(t)\|\leq \|e^{A_0 t}\| \|E(0)\| $. In order to obtain an upper bound for $ \|e^{A_0 t}\| $, a Lyapunov  function candidate $ L=\Tr[e^{A_0^T t}P_0e^{A_0 t}] $ with $ P_0 A_0+A_0^T P_0=-I $ is considered. Hence it can be obtained that $ \dot L=-\Tr[e^{A_0^T t}e^{A_0 t}] $. On the other hand, we have
\begin{equation}\label{L-ineq}
\underline{\lambda}(P_0)\Tr[e^{A_0^T t}e^{A_0 t}]\leq L(t)\leq \overline{\lambda}(P_0)\Tr[e^{A_0^T t}e^{A_0 t}]
\end{equation}
where $ \underline{\lambda}(P_0) $ and $ \overline{\lambda}(P_0) $ are the smallest and largest eigenvalues of $ P_0 $, respectively. One can use the preceding equation to get $ \dot L\leq-({1}/{\overline{\lambda}(P_0)})L $. Now, it is valid to say
\begin{equation*}
L(t)\leq e^{-\frac{1}{\overline{\lambda}(P_0)}t}L(0)
\end{equation*}
Moreover, by using \eqref{L-ineq}, one has
\begin{equation*}
\Tr[e^{A_0^T t}e^{A_0 t}]\leq n\frac{\overline{\lambda}(P_0)}{\underline{\lambda}(P_0)} e^{-\frac{1}{\overline{\lambda}(P_0)}t}
\end{equation*}
In addition, the equation $ \|e^{A_0 t}\|^2\leq \Tr[e^{A_0^T t}e^{A_0 t}] $ can be employed for rewriting the preceding equation as follows
\begin{equation}\label{ineq-tm}
\|e^{A_0 t}\|\leq k e^{-\lambda t}
\end{equation}
where $ k=\sqrt{n{\overline{\lambda}(P_0)}/{\underline{\lambda}(P_0)}} $ and $ \lambda=1/(2\overline{\lambda}(P_0)) $. By using the obtained upper bound of $ \|e^{A_0 t}\| $, we get
\begin{equation}\label{ineq-E}
\|E(t)\|\leq k e^{-\lambda t} \|E(0)\|
\end{equation}
For $ R(t) $, the fact that it is a positive definite matrix can be employed. On the other hand, from \eqref{RLS} it can be seen that $ \dot{R}(t)\leq 0 $; therefore, by using $ R(0)=\mu I $, we have $ 0\leq R(t)\leq \mu I $, i.e., $ R(t) $ is bounded. Now it can be said that
\begin{equation}\label{ineq-R}
\|R(t)\| \leq \mu
\end{equation}
For obtaining the preceding equation, $ \|R(t)\|^2=\overline{\lambda}(R(t))^2 $ with  the largest eigenvalue of $ R(t) $ as $\overline{\lambda}(R(t)) $ is employed, which is valid since $ R(t) $ is a symmetric matrix.

Now, the obtained upper bounds can be used for proving the stability of the observation error system. In this regard, the  state estimation  \eqref{last_msmo} and the definition of $ E $ can be utilized to obtain 
\begin{equation}\label{reg}
\tilde x_o=E\hat{\bar\alpha} +  \tilde{x}_N(\hat{\bar\alpha})
\end{equation}
where $ \tilde x_o=\hat x_o-x $ and $ \tilde x_N=\hat x_N-x $. By using the previous equation and \eqref{RLS}, we have
\begin{equation}\label{alpha-dot}
\dot{\hat{\bar\alpha}}=-R E^TC^TC \tilde x_o
\end{equation}
Moreover, the observation error system of the $ N $th observer can be obtained using \eqref{sys} and \eqref{last_msmo} as follows
\begin{equation}\label{n-th-error}
\dot{\tilde{x}}_N(\hat{\bar\alpha})=A_0\tilde{x}_N(\hat{\bar\alpha})+G_n\nu(\hat{x}_o)-D\xi(t,x,u)
\end{equation}
By taking the derivative of \eqref{reg} and using \eqref{e-dot}, \eqref{alpha-dot}, and \eqref{n-th-error}, one can obtain the observation error system as follows
\begin{equation}\label{reg-dot}
\dot{\tilde x}_o=A_0 \tilde x_o+G_n\nu(\hat{x}_o)-D\xi(t,x,u)-ER E^TC^TC \tilde x_o
\end{equation}
To show that the   error system \eqref{reg-dot} is stable, let us  start with the   the following nominal error system 
\begin{equation*}
\dot{\tilde x}_o=A_0 \tilde x_o+G_n\nu(\hat{x}_o)-D\xi(t,x,u)
\end{equation*}
One can compare the preceding equation with the observation error system for the SMO \eqref{smo}  and see that they are the same. As a result, by choosing the design matrices as \eqref{gl-gn}, the error system is quadratically stable \cite{edwards1994development}. Hence, there exists a Lyapunov function $ V(\tilde x_o) $ for the nominal error systems satisfying   
\begin{equation}
\begin{split}
c_1\|\tilde{x}_o\|^2 \leq V(\tilde x_o)&\leq c_2\|\tilde{x}_o\|^2\\
\frac{\partial V}{\partial t}+\frac{\partial V}{\partial \tilde{x}_o}[A_0 \tilde x_o&+G_n\nu(\hat{x}_o)-D\xi(t,x,u)]\leq -c_3\|\tilde{x}_o\|^2 \\
\|\frac{\partial V}{\partial \tilde{x}_o}\|\leq c_4\|\tilde{x}_o\|
\end{split} 
\label{lyap-nom1}
\end{equation}
for some positive constants $ c_1 $, $ c_2 $, $ c_3 $, and $ c_4 $ \cite{khalil1996nonlinear}. By considering $ V(\tilde x_o)$ as a Lyapunov function candidate for the observation error \eqref{reg-dot}, one can get

\begin{equation*}
\begin{split}
\dot V(\tilde x_o)&\leq \frac{\partial V}{\partial t}+\frac{\partial V}{\partial \tilde{x}_o}[A_0 \tilde x_o +G_n\nu(\hat{x}_o)-D\xi(t,x,u)]\\
&- \frac{\partial V}{\partial \tilde x_o} ERE^TC^TC\tilde{x}_o       
\end{split} 
\end{equation*}
Using    \eqref{lyap-nom1} and performing some basic manipulations on the preceding equation  result in 

\begin{equation*}
\dot V(\tilde x_o)\leq -c_3\|\tilde{x}_o\|^2+c_4 \|E\|^2 \|R\| \|C\|^2 \|\tilde{x}_o\|^2
\end{equation*}
One can employ \eqref{ineq-E} and \eqref{ineq-R} to get
\begin{equation*}
\dot V(\tilde x_o)\leq -c_3\|\tilde{x}_o\|^2+c_5 e^{-2\lambda t} \|\tilde{x}_o\|^2
\end{equation*}
where $ c_5=c_4 k^2\mu  \|E(0)\|^2  \|C\|^2 $. By considering \eqref{lyap-nom1}, we have
\begin{equation*}
\dot V(\tilde x_o)\leq (-\frac{c_3}{c_2}+\frac{c_5}{c_1} e^{-2\lambda t}) V(\tilde x_o)
\end{equation*}
Therefore, the following equation is satisfied
\begin{equation*}
V(t)\leq  e^{-\frac{c_3}{c_2}t}e^{\frac{c_5}{2\lambda c_1}(1- e^{-2\lambda t})} V(0)
\end{equation*}
Moreover, $ e^{\frac{c_5}{2\lambda c_1}(1- e^{-2\lambda t})}\leq e^{\frac{c_5}{2\lambda c_1}} $; hence 
\begin{equation*}
V(t)\leq k_1 e^{-\frac{c_3}{c_2}t} V(0)
\end{equation*}
where $ k_1=e^{\frac{c_5}{2\lambda c_1}} $. Moreover, one can employ \eqref{lyap-nom1} and obtain
\begin{equation}\label{exp-er}
\|\tilde{x}_o(t)\|\leq k_2 e^{-\frac{c_3}{2c_2}t} \|\tilde{x}_o(0)\|
\end{equation}
with $ k_2=\sqrt{k_1\frac{c_2}{c_1}} $. From the previous equation it can be concluded that the observation error converges to zero exponentially fast.

For proving the boundedness of $ \hat{\bar \alpha} $, the integral of \eqref{alpha-dot} is considered as follows
\begin{equation*}
\hat{\bar\alpha}(t)=\hat{\bar\alpha}(0)-\int_{0}^{t}R(\tau) E(\tau)^TC^TC \tilde x_o(\tau)d\tau
\end{equation*}
Now, one can use \eqref{ineq-E}, \eqref{ineq-R}, and \eqref{exp-er} to get
\begin{equation*}
\begin{split}
\|\hat{\bar\alpha}(t)\| &\leq \|\hat{\bar\alpha}(0)\|+\mu k k_2\|E(0)\|\|C\|^2 \|\tilde{x}_o(0)\| \int_{0}^{t} e^{-\lambda \tau}  e^{-\frac{c_3}{2c_2}\tau} d\tau
\end{split}
\end{equation*}
By considering $ \int_{0}^{t} e^{-\lambda \tau}  e^{-\frac{c_3}{2c_2}\tau} d\tau\leq {2c_2}/(2\lambda c_2+c_3) $, it can be seen that $ \hat{\bar \alpha} $ is bounded.

For $ \hat{x}_i(\hat{\bar \alpha}) $s, one can consider the estimation error of the $ i $th observer as $\tilde{x}_i=\hat{x_i}-x  $ and a Lyapunov function candidate $ V_i(\tilde{x}_i(\hat{\bar \alpha}))=\tilde{x}_i(\hat{\bar \alpha})^TP_0\tilde{x}_i(\hat{\bar \alpha}) $ with $ P_0 A_0+A_0^T P_0=-I $. Therefore, by using \eqref{sys} and \eqref{last_msmo} we can get
\begin{equation*}
\begin{split}
\dot V_i(\tilde{x}_i(\hat{\bar \alpha}))&\leq - \|\tilde{x}_i(\hat{\bar \alpha})\|^2+2\|\tilde{x}_i(\hat{\bar \alpha})\|\|P_0\|\\
&\times(\|G_n\| \|\nu(\hat{x}_o)\|  +\|D\| \|\xi\|)
\end{split}
\end{equation*}
On the other hand, we have $ \|\nu(\hat{x}_o)\|\leq \rho $ and $ \|\xi\|\leq \rho $. As a result, $ \dot V_i(\tilde{x}_i(\hat{\bar \alpha}))<0 $ for $ \|\tilde{x}_i(\hat{\bar \alpha})\|>2\rho \|P_0\|(\|G_n\|  +\|D\| ) $. In other words, since $ x $ is bounded, $ \hat{x}_i $s are uniformly ultimately bounded. 
\hfil$ \square $

The presented theorem states that the proposed observation scheme is able to provide a state estimation that converges to the state of the plant. In the next section, the performance of this observation scheme is investigated to obtain conditions that result in a better state estimation.

\subsection{Performance Investigation}

This section is aimed at analyzing the performance of the proposed observer to see whether it is able to provide better state estimations than a single sliding mode observer. In order to investigate this improvement, the following lemma  needs  to be considered.

\begin{lemma}\label{lem:obser}
	For system \eqref{sys}, if the pair $ (A,C) $ is observable, then the system $ (A,D,C) $ has no invariant zero.
\end{lemma}

\textbf{Proof.} According to the Popov-Belevitch-Hautus rank test, the following matrix has full column rank since $ (A,C) $ is observable
\begin{equation*}
\mathbb{P}(s)=\begin{bmatrix}
sI_n-A\\
C
\end{bmatrix}
\end{equation*}
Therefore, $ \text{rank}(\mathbb{P}(s))=n $. On the other hand, by considering the Rosenbrock matrix \eqref{rosen} and the fact that $ D $ has full rank, one can conclude that $ \text{rank}(\mathbb{R}(s))=\text{rank}(\mathbb{P}(s))+q $, which means  there is no $ s_0 $ that makes $ \mathbb{R}(s_0) $ lose rank. \hfill $ \square $


It can be seen from Lemma~\ref{lem:obser} that if $ \text{rank}(CD)=q $ and $ (A,C) $ is observable, the conditions of Theorem~\ref{th:smo-cond} are satisfied, and in turn, the results of the previous sections are still valid. As a result, in this section, it is assumed that the pair $ (A,C) $ is observable. In addition, we assume  that  $ x(0) $ is in the convex hull $ \mathcal{K} $ of $ \hat{x}_i(0) $s; hence \eqref{reg-0} is satisfied, and one can rewrite \eqref{reg} as follows
\begin{equation}\label{reg-per}
\tilde x_o=E\tilde{\bar\alpha} +  \hat{x}_N(\hat{\bar\alpha})-\hat{x}_N({\bar\alpha}^*)
\end{equation}
where $ \tilde{\bar\alpha}=\hat{\bar\alpha}-{\bar\alpha}^* $. For analyzing the preceding equation, the definition of $ \nu $ and the fact that by choosing $ \alpha_i=\alpha_i^* $s, the equality $ \hat{x}_o(t)=x(t) $ is valid, can be considered to get 
\begin{equation*}
\dot{\hat{x}}_N({\bar\alpha}^*)=A\hat{x}_N({\bar\alpha}^*)+Bu-G_le_N({\bar\alpha}^*)
\end{equation*} 
Therefore, we have
\begin{equation*}
\frac{d}{dt}[\hat{x}_N(\hat{\bar\alpha})-\hat{x}_N({\bar\alpha}^*)]=A_0[\hat{x}_N(\hat{\bar\alpha})-\hat{x}_N({\bar\alpha}^*)]+G_n\nu(\hat x_o)
\end{equation*}
Since $ \hat{x}_N(\hat{\bar\alpha}(0),0)=\hat{x}_N({\bar\alpha}^*,0) $, the following equation can be obtained
\begin{equation*}
\hat{x}_N(\hat{\bar\alpha}(t),t)-\hat{x}_N({\bar\alpha}^*,t)=\int_{0}^{t}e^{A_0(t-\tau)}G_n\nu(\hat x_o(\tau))d\tau
\end{equation*}
Now, by using \eqref{ineq-tm} and $ \|\nu(\hat x_o)\| \leq \rho $, one has
\begin{equation}\label{norm-sig}
\|\hat{x}_N(\hat{\bar\alpha}(t),t)-\hat{x}_N({\bar\alpha}^*,t)\|\leq\frac{1}{\lambda}k\rho \|G_n\|(1-e^{-\lambda t})
\end{equation}
The preceding equation and \eqref{reg-per} can be employed to get
\begin{equation}\label{reg-per1}
\|\tilde x_o(t)\|\leq\|E(t)\tilde{\bar\alpha}(t)\| +  \frac{1}{\lambda}k\rho \|G_n\|(1-e^{-\lambda t})
\end{equation}
To proceed with the analysis, it is required to consider $ \|E(t)\tilde{\bar\alpha}(t)\| $. In this regard, $ \eqref{alpha-dot} $ and \eqref{reg-per} can be utilized to get
\begin{equation}\label{alpha-tild-dot}
\dot{\tilde{\bar\alpha}}=-R E^TC^TC (E\tilde{\bar\alpha} +  \hat{x}_N(\hat{\bar\alpha})-\hat{x}_N({\bar\alpha}^*))
\end{equation}
On the other hand, one can use \eqref{RLS} and $ \frac{dR^{-1}}{dt}=-R^{-1}\dot{R}R^{-1} $ to obtain the following equation
\begin{equation}\label{r-inv-dot}
\frac{dR^{-1}}{dt}=E^T C^T CE
\end{equation}
Therefore, by employing \eqref{alpha-tild-dot} we have
\begin{equation*}
\frac{d(R^{-1}\tilde{\bar\alpha})}{dt}=-E^TC^TC(  \hat{x}_N(\hat{\bar\alpha})-\hat{x}_N({\bar\alpha}^*))
\end{equation*}
One can take the integral of the previous equation and premultiply it by $ R(t) $ to get
\begin{equation}\label{alpha-tild}
\tilde{\bar\alpha}=\frac{R(t)}{\mu}[\tilde{\bar\alpha}_0-\mu \int_{0}^{t}E(\tau)^TC^TC(  \hat{x}_N(\hat{\bar\alpha}(\tau),\tau)-\hat{x}_N({\bar\alpha}^*,\tau))d\tau]
\end{equation}
where $ \tilde{\bar\alpha}_0=\hat{\bar\alpha}_0-{\bar\alpha}^* $. It is worth noting that $ R(0)=\mu I $ is used for obtaining the preceding equation.

In what follows, the goal is to obtain $ \frac{1}{\mu}R(t) $ which exists in \eqref{alpha-tild}. Towards this end, \eqref{r-inv-dot} is employed to show that
\begin{equation*}
R(t)^{-1}-\frac{1}{\mu}I=\int_{0}^{t}E(\tau)^T C^T CE(\tau)d\tau
\end{equation*}
Using $ E(t)=e^{A_0t}E(0) $, we have
\begin{equation}\label{r-mu0}
\frac{R(t)}{\mu}=[I+\mu E(0)^T\int_{0}^{t} e^{A_0^T\tau}C^T Ce^{A_0\tau}d\tau E(0)]^{-1}
\end{equation}
Since it is assumed that $ (A,C) $ is observable, it can be easily shown that the pair $ (A_0,C) $ is also observable. As a result, the observability Gramian 
\begin{equation}\label{obs-gram}
W_o(t)=\int_{0}^{t} e^{A_0^T\tau}C^T Ce^{A_0\tau}d\tau
\end{equation}
is nonsingular for any $ t>0 $. It follows that there exists a constant $ b_1>0 $ such that $ b_1I\leq W_o(t) $. Hence $ W_o(t)^{-1} \leq \frac{1}{b_1}I$; and since $ W_o(t)=W_o(t)^T $, we have $ \|W_o(t)^{-1}\|^2=\overline\lambda (W_o(t)^{-1})^2$ which follows that
\begin{equation}\label{uper-obs-gram}
\|W_o(t)^{-1}\|\leq \frac{1}{b_1}
\end{equation}
Now by substituting \eqref{obs-gram} into \eqref{r-mu0} and using the matrix inversion lemma, one has
\begin{equation*}
\frac{1}{\mu}R(t)=I-\mu E(0)^T[W_o(t)^{-1}+\mu E(0)E(0)^T]^{-1} E(0)
\end{equation*}
By using the fact that $ E(0) $ is $ n\times(N-1) $ and $ N-1\geq n $, we can assume that $ E(0)E(0)^T $ is invertible. Therefore, we have
\begin{equation}\label{r-mu1}
\frac{1}{\mu}R(t)=I- E(0)^T(E(0)E(0)^T)^{-1}Q_0(t) E(0)
\end{equation}
where
\begin{equation}\label{q0}
Q_0(t)=[I+\frac{1}{\mu}W_o(t)^{-1}(E(0)E(0)^T)^{-1}]^{-1}
\end{equation}
It is required to rewrite \eqref{q0} into an infinite series. In this regard,  we will use Lemma \ref{lem:neum}.
In order to employ the   Lemma  \ref{lem:neum}, \eqref{uper-obs-gram} and $ \|(E(0)E(0)^T)^{-1}\|=b_2 $ can be used to get
\begin{equation*}
\|W_o(t)^{-1}(E(0)E(0)^T)^{-1}\|\leq \frac{b_2}{b_1}
\end{equation*}
Thus by choosing $ \mu>b_2/b_1 $, the condition of Lemma~\ref{lem:neum} is satisfied and \eqref{q0} can be considered as 
\begin{equation*}
Q_0(t)=I+\sum_{k=1}^{\infty}[-\frac{1}{\mu}W_o(t)^{-1}(E(0)E(0)^T)^{-1}]^{k}
\end{equation*}
By considering \eqref{reg-per1} and \eqref{alpha-tild}, for obtaining $ \|E(t)\tilde{\bar \alpha}(t)\| $, one needs to substitute the preceding equation into \eqref{r-mu1} and use $ E(t)=e^{A_0t}E(0) $ to get
\begin{equation*}
\frac{1}{\mu}E(t)R(t)=-e^{A_0t}\sum_{k=1}^{\infty}[-\frac{1}{\mu}W_o(t)^{-1}(E(0)E(0)^T)^{-1}]^{k}E(0)
\end{equation*}
As a result, using \eqref{alpha-tild} we have
\begin{equation}
\begin{split}
E(t)\tilde{\bar \alpha}(t)&=Q_1(t)E(0)\tilde{\bar\alpha}_0
\\
&-Q_2(t)  \int_{0}^{t}e^{A_0^T\tau}C^TC(  \hat{x}_N(\hat{\bar\alpha}(\tau),\tau)-\hat{x}_N({\bar\alpha}^*,\tau))d\tau
\end{split}\label{e-alpha}
\end{equation}
where 
\begin{align*}
Q_1(t)&=-e^{A_0t}\sum_{k=1}^{\infty}[-\frac{1}{\mu}W_o(t)^{-1}(E(0)E(0)^T)^{-1}]^{k}\\
Q_2(t)&=e^{A_0t}W_o(t)^{-1}(I+\sum_{k=1}^{\infty}[-\frac{1}{\mu}(E(0)E(0)^T)^{-1}W_o(t)^{-1}]^{k})
\end{align*}
Using $ \mu>b_2/b_1 $ and \eqref{ineq-tm} follows that
\begin{align*}
\|Q_1(t)\|&\leq ke^{-\lambda t}\sum_{k=1}^{\infty}[\frac{1}{\mu}\frac{b_2}{b_1}]^k=k\frac{b_2}{\mu b_1-b_2}e^{-\lambda t}\\
\|Q_2(t)\|&\leq ke^{-\lambda t}\frac{1}{b_1}(1+\sum_{k=1}^{\infty}[\frac{1}{\mu}\frac{b_2}{b_1}]^k)=k\frac{\mu}{\mu b_1-b_2}e^{-\lambda t}
\end{align*}
Utilizing the obtained upper bounds, \eqref{ineq-tm}, and \eqref{e-alpha} results
\begin{equation*}
\begin{split}
\|E(t)\tilde{\bar \alpha}(t)\|&\leq k\frac{b_2}{\mu b_1-b_2}\|E(0)\tilde{\bar\alpha}_0\|e^{-\lambda t} \\
&+k^2\frac{\mu}{\mu b_1-b_2} \|C\|^2 e^{-\lambda t} \\
&\times\int_{0}^{t}e^{-\lambda\tau}\|(  \hat{x}_N(\hat{\bar\alpha}(\tau),\tau)-\hat{x}_N({\bar\alpha}^*,\tau))\|d\tau
\end{split}
\end{equation*}
Now, \eqref{norm-sig} and \eqref{reg-per1} are employed to get
\begin{equation*}
\begin{split}
\|\tilde x_o(t)\|&\leq k\frac{b_2}{\mu b_1-b_2}\|E(0)\tilde{\bar\alpha}_0\|e^{-\lambda t}  +\\& \frac{1}{\lambda^2} \frac{\mu}{\mu b_1-b_2} k^3\rho\|C\|^2 \|G_n\| e^{-\lambda t}\times \\
&(\frac{1}{2}-e^{-\lambda t}+\frac{1}{2}e^{-2\lambda t})+ \frac{1}{\lambda}k\rho \|G_n\|(1-e^{-\lambda t})
\end{split}
\end{equation*}
Finally, since $ 1-e^{-\lambda t}\leq 1 $ and $ e^{-\lambda t}({1}/{2}-e^{-\lambda t}+e^{-2\lambda t}/2)\leq 2/27 $, we have
\begin{equation}\label{reg-per2}
\begin{split}
\|\tilde x_o(t)\|&\leq k\frac{b_2}{\mu b_1-b_2}\|E(0)\tilde{\bar\alpha}_0\|
 + \frac{2}{27}\frac{1}{\lambda^2} \frac{\mu}{\mu b_1-b_2} k^3\rho\|C\|^2 \|G_n\| 
 \\& +\frac{1}{\lambda}k\rho \|G_n\|
\end{split}
\end{equation}
For performing the comparison, it is required to analyze the observation error of conventional SMO using a similar method. In this regard, let us to assume that $ \dot{\hat{\bar \alpha}}(t)=0 $. Hence we have $ {\hat{\bar \alpha}}(t)={\hat{\bar \alpha}}_0 $, and one can use \eqref{last_msmo} and consider a state estimation as follows
\begin{equation*}
\hat x_s(t)=\sum\limits_{i=1}^{N-1} \hat{\alpha}_i(0) \hat{x}_i(\hat{\bar\alpha}_0,t) + (1-\sum\limits_{i=1}^{N-1} \hat{\alpha}_i(0)) \hat{x}_N(\hat{\bar\alpha}_0,t)
\end{equation*}
By using \eqref{last_msmo} it can be concluded that
\begin{equation*}
\dot{\hat x}_s=A{\hat x}_s+Bu-G_le_{os}+G_n\nu({\hat x}_s)
\end{equation*}
where $ e_{os}=C{\hat x}_s-y $. By comparing the preceding equation with \eqref{smo}, it can be seen that if $ \hat x_s(0)=\hat x(0) $, the obtained estimation from conventional SMO is equal to $ \hat x_s(t) $. Therefore, by using \eqref{reg-per} and considering $ \dot{\hat{\bar \alpha}}(t)=0 $, one can get
\begin{equation*}
\tilde x(t)=E(t)\tilde{\bar\alpha}_0 +  \hat{x}_N(\hat{\bar\alpha}_0,t)-\hat{x}_N({\bar\alpha}^*,t)
\end{equation*}
Similar to \eqref{reg-per1} it can be obtained
\begin{equation*}
\|\tilde x(t)\|\leq \|E(t)\tilde{\bar\alpha}_0\| +  \frac{1}{\lambda}k\rho \|G_n\|(1-e^{-\lambda t})
\end{equation*}
Moreover, employing $ E(t)=e^{A_0t}E(0) $ and \eqref{ineq-tm} results
\begin{equation}\label{norm-singl}
\|\tilde x(t)\|\leq k \|E(0)\tilde{\bar\alpha}_0\| +  \frac{1}{\lambda}k\rho \|G_n\|
\end{equation}

Finally, it can be seen that the right hand side of \eqref{reg-per2} can become arbitrary smaller than the upper bound of $ \|\tilde{x}(t)\| $ in \eqref{norm-singl}  if we choose $ \mu $ and $ \lambda $ big enough. In other words, based on the definition of $ \lambda $, i.e., $ \lambda=1/(2\overline{\lambda}(P_0)) $ where  $ P_0 A_0+A_0^T P_0=-I $, one can conclude that choosing   $ |\text{Re}(\lambda_i(A_0))| $s and $ \mu $ big enough can result in a better estimation in comparison to conventional sliding mode observers.

\subsection{Fault Detection and Isolation}
In this section, the fault of system \eqref{sys} is approximated using the proposed observer. In this regard, it is required to show that a sliding motion takes place on the surface in the error space
\begin{equation}\label{S}
\mathcal S=\{\tilde{x}_o\in\mathcal{R}^n:C\tilde{x}_o=0\}
\end{equation}
Towards this end, the following lemma is presented.
\begin{lemma}\label{lem:sm}
	Let all conditions of Theorem~\ref{th:conver} are satisfied. Then, a sliding motion takes place on $ \mathcal S $ \eqref{S} in finite time.
\end{lemma}

\textbf{Proof.} By using the transition matrix $J$, one can transform the error system \eqref{reg-dot} to the following form
\begin{equation}
\begin{split}
\dot{\tilde x}_{I}&=\mathcal{A}_{11}{\tilde x}_{I}-\mathcal{C}_1^T JERE^T J^T \mathcal{C} {\tilde y}_{o} \\
\dot{\tilde y}_{o}&=\mathcal A_{21}{\tilde x}_{I}+\mathcal A_{22}^s {\tilde y}_{o}+\|\mathcal D_2\|\nu (\hat{x}_o)\\
&-\mathcal D_2 \xi -\mathcal{C}^T JERE^T J^T \mathcal{C} {\tilde y}_{o}
\end{split}
\label{er-sm1}
\end{equation}
where $ \mathcal{C}_1=\begin{bmatrix}
I_{n-p} & 0
\end{bmatrix}^T $ and $ \mathcal{C}=\begin{bmatrix}
0 & I_p
\end{bmatrix}^T $. 
Now, $ V_s({\tilde y}_{o})= {\tilde y}_{o}^T P_2  {\tilde y}_{o}$ with $ P_2 $ satisfying \eqref{lyp-eq} is considered to get
\begin{equation*}
\begin{split}
\dot V_s&=- {\tilde y}_{o}^T Q_2 {\tilde y}_{o}-2 {\tilde y}_{o}^T P_2\mathcal{C}^T JERE^T J^T \mathcal{C} {\tilde y}_{o}\\
&+2 {\tilde y}_{o}^T P_2 (\mathcal A_{21}{\tilde x}_{I}+\|\mathcal D_2\|\nu (\hat{x}_o)-\mathcal D_2 \xi )
\end{split}
\end{equation*}
One can use $ \|\mathcal{C}\|=1 $, \eqref{ineq-E}, and \eqref{ineq-R} to obtain the following  equation
\begin{equation*}
\begin{split}
\dot V_s&\leq (-\underline{\lambda}(Q_2)+2 \mu k^2 \|P_2\| \|J\|^2 \|E(0)\|^2 e^{-2\lambda t}) \|{\tilde y}_{o}\|^2\\
&+2 {\tilde y}_{o}^T P_2 (\mathcal A_{21}{\tilde x}_{I}+\|\mathcal D_2\|\nu (\hat{x}_o)-\mathcal D_2 \xi )
\end{split}
\end{equation*}
It can be seen that there exists $ T_0>0 $ such that for $ t\geq T_0 $ we have
\begin{equation*}
-\underline{\lambda}(Q_2)+2 \mu k^2 \|P_2\| \|J\|^2 \|E(0)\|^2 e^{-2\lambda t}\leq 0
\end{equation*}
Hence, it is valid to say that
\begin{equation*}
\dot V_s\leq 2 {\tilde y}_{o}^T P_2 (\mathcal A_{21}{\tilde x}_{I}+\|\mathcal D_2\|\nu (\hat{x}_o)-\mathcal D_2 \xi )\qquad  \forall T_0\leq t
\end{equation*}
From the definition of $ \nu $, it can be seen that ${\tilde y}_{o}^T P_2 \nu (\hat{x}_o)=-\rho \|P_2 {\tilde y}_{o}\|  $. Moreover, from \eqref{up-xi} and the definition of $ \rho $, we have $ \|\xi\|\leq \rho-\gamma_0 $. Therefore, it is obtained 
\begin{equation*}
\dot V_s\leq 2\|P_2 {\tilde y}_{o}\| (\|\mathcal A_{21}{\tilde x}_{I}\|-\|\mathcal D_2 \| \gamma_0)\qquad  \forall T_0\leq t
\end{equation*}
It can be seen from the preceding equation that in the domain $ \Omega=\{{\tilde x}_{I}: \|\mathcal A_{21}{\tilde x}_{I}\|<\|\mathcal D_2 \| \gamma_0-\delta\} $ with  positive scalar $ \delta $, the following equation is valid
\begin{equation*}
\dot V_s\leq -2\delta \|P_2 {\tilde y}_{o}\|
\end{equation*}
From Theorem~\ref{th:conver}, we know that the observation error system is quadratically stable; hence  $ {\tilde x}_{I} $ enters $ \Omega $ in finite time. As a result, there exists $ T_0\leq T_1 $ such that for all $ t\geq T_1 $ the preceding equation is valid. In addition, since $ \|P_2 {\tilde y}_{o}\|^2= (P_2^{1/2} {\tilde y}_{o})^T P_2(P_2^{1/2} {\tilde y}_{o})$, it can obtained
\begin{equation*}
\underline{\lambda}(P_2) V_s\leq \|P_2 {\tilde y}_{o}\|^2
\end{equation*}
By employing the previous equation, one can get
\begin{equation*}
\dot V_s\leq -2\delta \sqrt{\underline{\lambda}(P_2)} \sqrt{V_s} \qquad  \forall T_1\leq t
\end{equation*}
Solving the obtained inequality results
\begin{equation*}
\sqrt{V_s(t)}\leq \sqrt{V_s(T_1)}-\delta \sqrt{\underline{\lambda}(P_2)} (t-T_1)
\end{equation*}
It can be easily seen that if $ T_s=T_1+\sqrt{V_s(T_1)/(\delta \underline{\lambda}(P_2))} $, $V_s(t)=0  $ for all $ t\geq T_s $. In another word, sliding motion takes place on $ \mathcal{S} $ in finite time.  \hfil$ \square $

The presented lemma can be employed for fault reconstruction of \eqref{sys}. 
Since sliding motion takes place after a finite time, we have $ {\tilde y}_{o}=0 $ and $ \dot{\tilde y}_{o}=0 $, and  equation \eqref{er-sm1}  become
\begin{align}
\dot{\tilde x}_{I}&=\mathcal{A}_{11}{\tilde x}_{I} \label{er-sm1-ts}\\
0&=\mathcal A_{21}{\tilde x}_{I}+\|\mathcal D_2\|\nu_{{eq}} (\hat{x}_o)-\mathcal D_2 \xi \label{er-sm2-ts}
\end{align}
where $ \nu_{{eq}} $ is the \textit{equivalent signal}. According to the equivalent control method,  the equivalent signal helps the sliding motion to be maintained \cite{utkin2013sliding}. Moreover, since $ \mathcal{A}_{11} $ is Hurwitz, it can be seen from \eqref{er-sm1-ts} that $ {\tilde x}_{I}(t)\to 0 $ as $ t\to \infty $. Consequently, from \eqref{er-sm2-ts} one can get
\begin{equation*}
\|\mathcal D_2\|\nu_{{eq}} (\hat{x}_o)\to \mathcal D_2 \xi 
\end{equation*}
It is worth noting that $ \text{rank}(\mathcal D_2)=q $ \cite{edwards1994development}. Therefore, we can reconstruct  the fault as follows
\begin{equation}\label{f-hat}
\hat \xi_o =\|\mathcal D_2\| \mathcal D_2(\mathcal D_2^T\mathcal D_2)^{-1}\mathcal D_2^T \nu_{{eq}} (\hat{x}_o)
\end{equation}
For obtaining the equivalent signal, the presented method in \cite{edwards2000sliding} can be considered. In this regard, the discontinuous term $ \nu $ is replaced by its continuous approximation as follows
\begin{equation}\label{nu-d}
\nu_{\delta}(\hat{x}_o)=-\rho\frac{P_2 \tilde{y}_o}{\|P_2 \tilde{y}_o\|+\delta}
\end{equation}
where $ \delta $ is a positive constant that determines the accuracy of approximation and needs to be chosen sufficiently small. Then, the equivalent signal can be obtained to any desired accuracy by using $ \nu_{{eq}}=\nu_{\delta} $.

It is worth noting that by considering \eqref{f-hat}, one can see that the obtained fault reconstruction is based on the state estimation. On the other hand, it was shown in the preceding section that the proposed scheme is able to provide a more preferable state estimation.  Therefore, one can conclude that  it can also result in a better fault reconstruction in comparison to conventional sliding mode observers.
%
%

\section{Simulation Results}
In this section, the effectiveness of the proposed fault detection and reconstruction  methodology is illustrated through simulation. Towards this end, consider the Matsumoto-Chua-Kobayashi (MCK) circuit, which is a chaotic system, as follows \cite{tamasevicius1997hyperchaotic}:
\begin{equation}\label{mck}
\begin{split}
\dot{x}&=Ax+D\xi \\
y&=Cx
\end{split}
\end{equation}
where
\begin{equation}\label{mck-matrix}
\begin{split}
\xi(x)&=\begin{cases}
0.2+3(x_1-x_3+1) & x_1-x_3<-1\\
-0.2(x_1-x_3) & -1\leq x_1-x_3 \leq 1\\
-0.2+3(x_1-x_3-1) & 1<x_1-x_3
\end{cases}\\
A&=\begin{bmatrix}
0 & -1 & 0 & 0\\
1 & 0.7 & 0 & 0\\
0 & 0 & 0 & -10\\
0 & 0 & 1.5 & 0
\end{bmatrix},B=0,D=\begin{bmatrix}
-1\\0\\10\\0
\end{bmatrix}\\
C&=\begin{bmatrix}
1 & 0 & 0 & 0\\
0 & 1 & 0 & 1
\end{bmatrix}
\end{split}
\end{equation} 
The system initial condition is considered as $ x(0)=\begin{bmatrix}
-0.1 & 0 & 0.2 & 0
\end{bmatrix}^T $. 

In order to estimate the state variables and  reconstruct the fault rapidly, it is required to  design the matrices $ G_l $ and $ G_n $ in the proposed methodology. Towards this end, the eigenvalues of $ \mathcal{A}_{11} $ are selected as $ \{-4,-6\} $, and the presented algorithm in  \cite{edwards1994development} is utilized to obtain the following transformation matrix
\begin{equation*}
J=\begin{bmatrix}
10 &	22.857 &	1 &	22.857\\
0 &	-63.265 &	0 &	-64.265\\
1 &	0 &	0 &	0\\
0 &	1 &	0 &	1
\end{bmatrix}
\end{equation*}
Now, one can employ the aforementioned transition matrix  and transform \eqref{mck} into the equivalent form  \eqref{sys-trans}. Consequently, matrices $ (A,D,C) $ defined in \eqref{mck-matrix} are transformed into new coordinates $ (\mathcal{A},\mathcal{D},\mathcal{C}) $ as follows
\begin{equation*}
\begin{split}
\mathcal{A}&=\left[\begin{array}{c c | c c}
34.285 &  16 & -320 &  234.571\\
-96.398 &  -44.285 &  900.714 & -642.653\\
\cline{1-4}
0 &  -1 &         0 & -64.265 \\
1.5 &    0.7 &  -14 &   10.7 
\end{array}\right],\\
\mathcal{D}&=\left[\begin{array}{c}
0 \\ 0 \\ \hline -1 \\ 0
\end{array}\right],
\mathcal{C}=\left[\begin{array}{c c | c c}
0 & 0 & 1 & 0\\
0 & 0 & 0 & 1
\end{array}\right]
\end{split}
\end{equation*}
where $ \mathcal{A}=JAJ^{-1},\mathcal{D}=JD,\mathcal{C}=CJ^{-1} $.   Then by considering $ \mathcal{A}_{22}^s=-10I_2 $ and employing \eqref{smo-trans} and \eqref{gain-part}, one can get
\begin{equation*}
\mathcal{G}_l=\begin{bmatrix}
-320 &  234.571\\
900.714 &-642.653\\
10  &-64.265\\
-14 &  20.7
\end{bmatrix},\mathcal{G}_n=\begin{bmatrix}
0 & 0\\
0 & 0\\
1 & 0\\
0 & 1
\end{bmatrix}
\end{equation*}
Finally, the design matrices can be obtained using $ G_l=J^{-1}\mathcal{G}_l $ and $ G_n=J^{-1}\mathcal{G}_n $. Moreover, $ P_2=I_2 $, $ \rho=10 $, and $ \delta=0.01 $ are employed to  design  $ \nu_{\delta} $ as \eqref{nu-d}.

As stated through the paper, the initial conditions $ \hat{x}_i(0) $s should be chosen such that $ x(0) $ lies in their convex hull $ \mathcal{K} $, and for performance improvement, $ E(0)E(0)^T $ needs to be invertible.   In this regard, five  sliding mode observers \eqref{last_msmo} ($ N=5 $) with the following initial conditions are considered
\begin{equation*}
\begin{split}
\hat{x}_1(0)&=\begin{bmatrix}
+1 & -1 & +1 & -1	\end{bmatrix}^T\\
\hat{x}_2(0)&=\begin{bmatrix}	-1 & +1 & -1 & +1
\end{bmatrix}^T\\
\hat{x}_3(0)&=\begin{bmatrix}
+1 & +1 & +1 & -1
\end{bmatrix}^T\\
\hat{x}_4(0)&=\begin{bmatrix}
+1 & -1 & -1 & +1
\end{bmatrix}^T\\
\hat{x}_5(0)&=\begin{bmatrix}
+1 & -1 & +1 & +1
\end{bmatrix}^T
\end{split}
\end{equation*}
Then, the   initial conditions of  the RLS algorithm \eqref{RLS} are chosen as $ \hat{\bar\alpha}_0=\begin{bmatrix}
0.2 & 0.2 & 0.2 & 0.2
\end{bmatrix}^T$ and $ R(0)=\mu I_4 $ with $ \mu=10^2 $.

In addition to the proposed scheme, one sliding mode observer with the same design matrices,  parameters, and initial conditions, i.e., $ \hat{x}(0)=\sum_{N=1}^{5} \hat{\alpha}_i(0) \hat{x}_i(0)=\begin{bmatrix}
0.6 & -0.2 & 0.2 & 0.2
\end{bmatrix}^T$, is utilized to obtain state estimation and fault reconstruction. This assists us to demonstrate that the proposed observation scheme results in  more accurate  estimations in comparison to single sliding mode observers. 

\begin{figure}
	\centering 
	\includegraphics	[width=3.25in, bb= 31 6 383 296, clip=true]{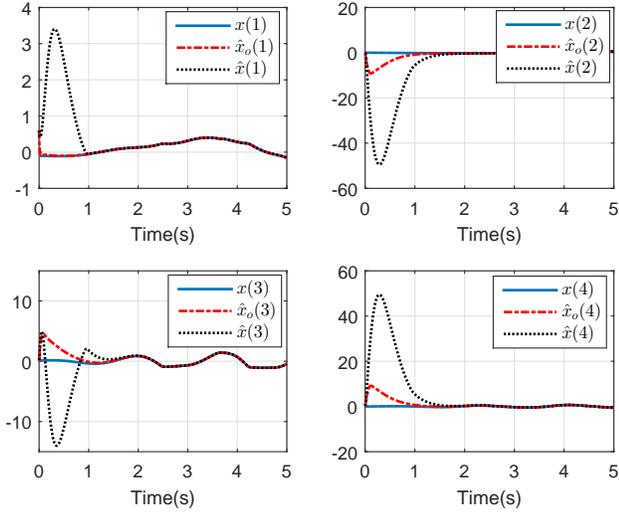} 
	\caption{State variables, $ x $, and their estimations using the proposed scheme, $ \hat{x}_o $, and single sliding mode observer, $ \hat{x} $, for $ \mu=10^{2} $.} 
	\label{fig:x_single_mu1}
\end{figure}

\begin{figure}
	\centering 
	\includegraphics
	[width=2.5in, bb= 39 1 383 296, clip=true]{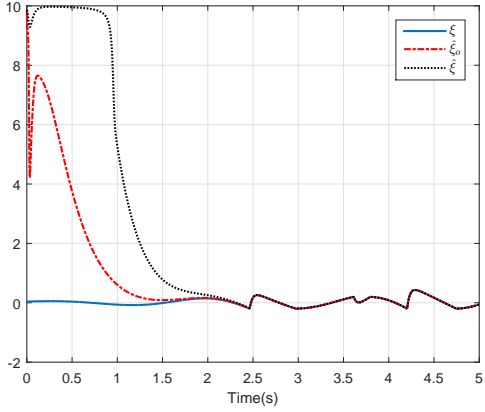} 
	\caption{Fault signal, $ \xi $, and its estimations using the proposed scheme, $ \hat{\xi}_o $, and single sliding mode observer, $ \hat{\xi} $, for $ \mu=10^{2} $.} 
	\label{fig:xi_single_mu1}
\end{figure}

The obtained simulation results are presented in Fig.~\ref{fig:x_single_mu1} and Fig.~\ref{fig:xi_single_mu1}.  From these figures, it can be easily seen that the consequence of employing the proposed observer is a state estimation with more preferable transient response.

%

\begin{figure}
	\centering 
	\includegraphics
	[width=3.25in, bb= 31 6 383 296, clip=true]{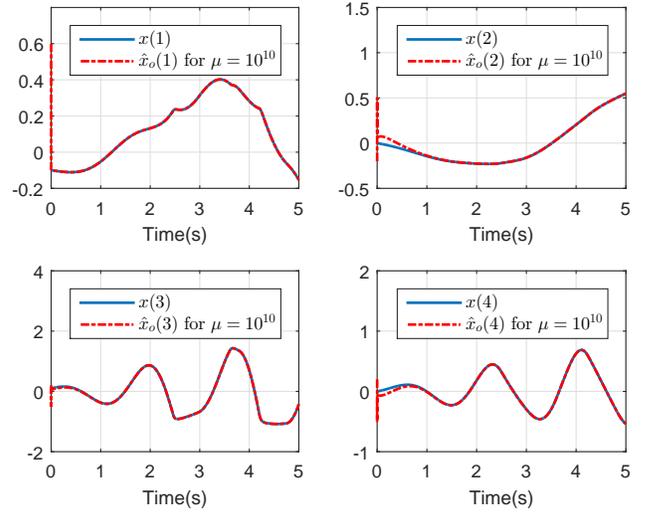} 
	\caption{State variables, $ x $, and their estimations using the proposed scheme, $ \hat{x}_o $, for $ \mu=10^{10} $.} 
	\label{fig:x_mu2}
\end{figure}

\begin{figure}
	\centering 
	\includegraphics
	[width=2.5in, bb= 39 1 383 296, clip=true]{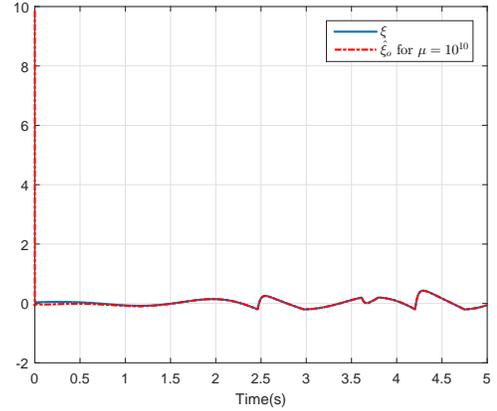} 
	\caption{Fault signal, $ \xi $, and its estimation using the proposed scheme, $ \hat{\xi}_o $, for $ \mu=10^{10} $.} 
	\label{fig:xi_mu2}
\end{figure}

In order to demonstrate that choosing $ \mu $ big enough results in a superior performance, the simulation is also performed for $ \mu=10^{10} $. Fig.~\ref{fig:x_mu2} and Fig.~\ref{fig:xi_mu2} show the obtained state and fault estimations for this value. The figures validate the theory and show that choosing $ \mu $ results in a better performance in comparison to the conventional SMO.

\section{Conclusions}
In this research a novel fault reconstruction  scheme was presented using the concept of second level adaptation   and SMOs. In the proposed scheme, the information provided by  multiple SMOs with suitably chosen initial conditions was  employed  to reconstruct  the system states  and the fault (unknown input) rapidly. In this regard, it was shown that if the initial condition of system lies inside the convex hull of SMOs  initial conditions, there exist some constant coefficients  that provide a perfect state estimation. An estimation of these  coefficients was obtained using the RLS algorithm. Mathematical analyses/justifications were provided to highlight performance of the proposed observation strategy.  The  stability of the overall system as well as the structure of the fault reconstruction scheme were fully addressed. Since the proposed approach employs the collective  information obtained from multiple SMOs, it results in  estimations with more preferable transient response in comparison to conventional SMO-based fault detection strategies. 

\end{document}